\documentclass[10pt]{article}

\usepackage{amsfonts}
\usepackage{amssymb}
\usepackage{amsmath}
\usepackage{enumerate}

\newtheorem{Lem}{\bf LEMMA}[section]

\newtheorem{Pro}{\bf PROPOSITION}[section]

\newtheorem{Cor}{\bf COROLLARY}[section]

\newfont{\grandsy}{cmsy10}

\title{THE NONSTATIONARY IDEAL ON $P_\kappa(\lambda)$ FOR $\lambda$ SINGULAR}

\author{Pierre MATET   \footnote{Publication 33.} \,
 and Saharon SHELAH \footnote{Research supported by the United States - Israel Binational Science Foundation (Grant no. 2002323). Publication 869.}}

\date{}

\begin{document}

\maketitle





\vskip 1,5cm

\begin{abstract} Let $\kappa$ be a regular uncountable cardinal and $\lambda>\kappa$ a singular
strong limit cardinal. We give a new characterization of the nonstationary subsets of
$P_\kappa(\lambda)$ and use this to prove that the nonstationary ideal on $P_\kappa(\lambda)$ is
nowhere precipitous.
\end{abstract}

\setcounter{section}{-1}

\section{Introduction}

Let $\kappa$ be a regular uncountable cardinal and 
$\lambda>\kappa$  a singular cardinal. Let $I_{\kappa,\lambda}$ (respectively,
$NS_{\kappa,\lambda})$ denote the ideal of noncofinal (respectively, nonstationary) subsets of
$P_\kappa(\lambda).$ Now suppose $\lambda$ is a strong limit cardinal. If $cf(\lambda)<\kappa,$
then by a result of Shelah [7], $NS_{\kappa,\lambda} = I_{\kappa,\lambda}\mid A$ for some $A.$ If
$cf(\lambda)\geq\kappa,$ then by results of [4], $NS_{\kappa,\lambda}\not= I_{\kappa,\lambda}\mid
A$ for every
$A.$ Nevertheless, Shelah's result can be generalized as follows. Given an infinite cardinal
$\mu\leq\lambda,$ let $NS_{\kappa,\lambda}^\mu$ denote the smallest $\mu$-normal ideal on
$P_\kappa(\lambda),$ where an ideal $J$ on $P_\kappa(\lambda)$ is said to be $\mu$-normal if for
every $A\in J^+$ and every $f:A\rightarrow\mu$ with the property that $f(a)\in a$ for all $a\in A$, there
exists $B\in J^+\cap P(A)$ with $f$ being constant on $B.$ Note that $NS_{\kappa,\lambda}^\lambda
= NS_{\kappa,\lambda},$ and $NS_{\kappa,\lambda}^\mu=I_{\kappa,\lambda}$ whenever $\mu<\kappa.$ We
will show that $NS_{\kappa,\lambda}= NS_{\kappa,\lambda}^{cf(\lambda)}\mid A$ for some $A.$ Since,
by a result of Matsubara and Shioyia [6], $I_{\kappa,\lambda}$ is nowhere precipitous, it
immediately follows that $NS_{\kappa,\lambda}$ is nowhere precipitous in case
$cf(\lambda)<\kappa,$ a result that is also due to Matsubara and Shioyia [6]. It is claimed in
[5] that $NS_{\kappa,\lambda}$ is also nowhere precipitous in case $cf(\lambda)\geq\kappa.$
Unfortunately, there is a mistake in the proof (see the last line of the proof of Lemma 2.9~:
$a\in C_{\alpha^*}[g_{\alpha^*}]$ does not necessarily imply that $a\in C[g]).$ We show that the
proof can be repaired by using our characterization of $NS_{\kappa,\lambda}.$

For the results above to hold, it is not necessary to assume that $\lambda$ is a strong limit
cardinal. In fact we show that if $\overline{\hbox{cof}}(NS_{\kappa,\tau})\leq\lambda$ for every
cardinal $\tau$ with $\kappa\leq\tau<\lambda,$ then
$NS_{\kappa,\lambda}=NS_{\kappa,\lambda}^{cf(\lambda)}\mid A$ for some $A.$ (If GCH holds
in $V$ and $\mathbb{P}$ is the forcing notion to add $\lambda^+$ Cohen reals, then in $V^\mathbb{P}, \lambda$ is
no longer a strong limit cardinal but, by results of [3], for every cardinal $\tau$ with
$\kappa\leq\tau<\lambda, \hbox{ cof}(NS_{\kappa,\tau})=\tau^+$ and hence $\overline{\hbox{
cof}}(NS_{\kappa,\tau}\leq\lambda).$ Let us observe that by results of [4], the converse holds in
case $cf(\lambda)<\kappa.$

Note that if $NS_{\kappa,\lambda} = NS_{\kappa,\lambda}^{cf(\lambda)}\mid A,$ then for each
cardinal $\chi$ with 	\ $\kappa\cdot(cf(\lambda))^+\leq\chi<\lambda, NS_{\kappa,\lambda}^\chi\mid A =
NS_{\kappa,\lambda}^{cf(\lambda)}|A.$ We show that for this (i.e. the existence of $A\in
NS_{\kappa,\lambda}^*$ such that $NS_{\kappa,\lambda}^\chi\mid A = 
NS_{\kappa,\lambda}^{cf(\lambda)}\mid A)$ to hold, it is sufficient to assume that
$\overline{\hbox{ cof}}(NS_{\kappa,\tau}^\chi)\leq\lambda$ for every cardinal $\tau$ with
$\chi\leq\tau<\lambda$.
\\

\section{Basic material}

{\bf Throughout the paper $\kappa$ denotes a regular  uncountable cardinal and $\lambda$ a
cardinal greater than or equal to $\kappa.$
}

$NS_\kappa$ denotes the nonstationary ideal on \ $\kappa.$

For a set $A$ and a cardinal $\rho,$ let $P_\rho(A) = \{a\subseteq A :\mid a\mid<\rho\}.$

$I_{\kappa,\lambda}$ denotes the set of all $A\subseteq P_\kappa(\lambda)$ such that $\{a\in
A:b\subseteq a\}=\phi$ for some $a\in P_\kappa(\lambda).$

By an {\sl ideal} on $P_\kappa(\lambda),$ we mean a collection $J$ of subsets of
$P_\kappa(\lambda)$ such that (i) $I_{\kappa,\lambda}\subseteq J$ ; (ii) $P_\kappa(\lambda)\notin
J$ ; (iii) $P(A)\subseteq J$ for all $A\in J$ ; and (iv) $\cup X\in J$ for every $X\in
P_\kappa(J).$

Given an ideal $J$ on $P_\kappa(\lambda),$ let $J^+ = \{A\subseteq P_\kappa(\lambda) : A\notin
J\}$ and $J^*=\{A\subseteq P_\kappa(\lambda) :P_\kappa(\lambda)\setminus A\in J\}.$ For $A\in
J^+,$ let $J\mid A = \{B\subseteq P_\kappa(\lambda) : B\cap A\in J\}.$ ${\cal M}_J$ denotes the
collection of all $Q\subseteq J^+$ such that  (i) $A\cap B\in J$ for any distinct $A, B\in J,$ and
(ii) for every $C\in J^+,$ there is $A\in Q$ with $A\cap C\in J^+.$ For a cardinal $\rho, J$ is
{\sl $\rho$-saturated} if $\mid Q\mid<\rho$ for every $Q\in{\cal M}_J.$

An ideal $J$ on $P_\kappa(\lambda)$ is {\sl precipitous} if whenever $A\in J^+$ and
$<Q_n:n<\omega>$ is a sequence of members of ${\cal M}_{J\mid A}$ such that
$Q_{n+1}\subseteq\displaystyle\bigcup_{B\in Q_n} P(B)$ for all $n<\omega,$ there exists $f\in
\displaystyle\prod_{n\in\omega}Q_n$ such that $f(0)\supseteq f(1)\supseteq\ldots$ and
$\displaystyle\bigcap_{n<\omega}f(n)\not=\phi.$ $J$ is {\sl nowhere precipitous} if for each $A\in
J^+, J\mid A$ is not precipitous. $G(J)$ denotes the following two-player game lasting $\omega$
moves, with player I making the first move : I and II alternately pick members of $J^+,$ thus
building a sequence $<X_n:n<\omega>,$ subject to the condition that $X_0\supseteq
X_1\supseteq\ldots$ II wins
$G(J)$ just in case $\displaystyle\bigcap_{n<\omega}X_n=\phi.$
\\
\begin{Lem}([2]) {\sl An ideal $J$ on $P_\kappa(\lambda)$ is nowhere precipitous if and only
if II has a winning strategy in the game $G(J).$
}
\end{Lem}

Given an ideal $J$ on $P_\kappa(\lambda), \hbox{ cof}(J)$ denotes the least cardinality of any
$X\subseteq J$ such that $J=\displaystyle\bigcup_{A\in X}P(A).$ $\overline{\hbox{ cof}}(J)$
denotes the least size of any $Y\subseteq J$ with the property that for every $A\in J,$ there is
$y\in P_\kappa(Y)$ with $A\subseteq\cup y.$ Let $u(\kappa,\lambda) = \hbox{
cof}(I_{\kappa,\lambda}).$ The following is well-known (see e.g. [3]) :
\\

\begin{Lem} $\lambda^{<\kappa}=2^{<\kappa}\cdot u(\kappa,\lambda).$
\end{Lem}

\begin{Lem}([3]){\sl Let $A\in I_{\kappa,\lambda}^+$ be such that $\mid\{a\in A :
b\subseteq a\}| = |A|$ for every $b\in P_\kappa(\lambda).$ Then $A$ can be decomposed into $|A|$
pairwise disjoint members of $I_{\kappa,\lambda}^+.$
}
It follows that if $NS_{\kappa,\lambda}=I_{\kappa,\lambda}\mid A$ for some $A,$ then (a)
$P_\kappa(\lambda)$ can be split into $c(\kappa,\lambda)$ disjoint stationary sets, where
$c(\kappa,\lambda)$ denotes the least size of any closed unbounded subset of $P_\kappa(\lambda),$
and (b) every stationary subset of $P_\kappa(\lambda)$ can be split into $u(\kappa,\lambda)$
disjoint stationary sets.
Let $\mu$ and $\theta$ be two cardinals such that $1\leq\mu\leq\lambda$ and
$2\leq\theta\leq\kappa.$ An ideal $J$ on $P_\kappa(\lambda)$ is  {\sl $[\mu]^{<\theta}$-normal} if
given $A\in J^+$ and $f : A\rightarrow P_\theta(\mu)$ with the property that $f(a)\in
P_{|a\cap\theta|}(a\cap\mu)$ for all $a\in A,$ there exists $B\in J^+\cap P(A)$ such that $f$ is
constant on $B$. (Note that $[\lambda]^{<\kappa}$-normality is the same as the well-known notion
of strong normality). {\bf We set $\overline\theta=\theta$ if $\theta<\kappa,$ or $\theta=\kappa$
and $\kappa$ is a limit cardinal, and $\overline\theta=\nu$ if $\theta=\kappa=\nu^+.$}
\end{Lem}

\begin{Lem} ([3]) 
\begin{enumerate}[i)]
\item   {\sl Suppose that $\mu<\kappa,$ or $\theta<\kappa,$ or $\kappa$ is not a limit cardinal.
Then there exists a $[\mu]^{<\theta}$-normal ideal on $P_\kappa(\lambda)$ if and only if
$|P_{\overline\theta}(\rho)|<\kappa$ for every cardinal $\rho<\kappa\cap(\mu+1).$
}

\item {\sl Suppose that $\mu\geq\kappa, \theta=\kappa$ and $\kappa$ is a limit cardinal. Then
there exists a $[\mu]^{<\theta}$-normal ideal on $P_\kappa(\lambda)$ if and only if $\kappa$ is a
Mahlo cardinal.
}
\end{enumerate}
\end{Lem}

Assuming there exists a $[\mu]^{<\theta}$-normal ideal on $P_\kappa(\lambda),
NS_{\kappa,\lambda}^{[\mu]^{<\theta}}$ denotes the smallest such ideal.

\begin{Lem}([3])
\begin{enumerate}[i)]
\item  {\sl Suppose $\mu<\kappa.$ Then
$NS_{\kappa,\lambda}^{[\mu]^{<\theta}}=I_{\kappa,\lambda}.$
}
\item  {\sl Suppose $\theta\leq\omega.$ Then
$NS_{\kappa,\lambda}^{[\mu]^{<\theta}}=NS_{\kappa,\lambda}^\mu.$
}
\end{enumerate}
\end{Lem}

For $g : P_{\overline\theta\cdot 3}(\mu)\rightarrow P_3(\lambda),$ let $C_g^{\kappa,\lambda}$ be the set
of all $a\in P_\kappa(\lambda)$ such that $a\cap(\overline\theta\cdot 3)\not=\phi$ and
$f(e)\subseteq a$ for every $e\in P_{|a\cap(\overline\theta\cdot 3)|}(a\cap\mu).$
\\  

\begin{Lem}([4]) {\sl Suppose $\kappa\leq\mu<\lambda<\mu^{+\kappa}.$ \ Then
$\overline{cof}(NS_{\kappa,\lambda}^\mu)=\lambda\cdot\overline{cof}(NS_{\kappa,\mu}).$
}
\end{Lem}

\begin{Lem}([3]) {\sl Suppose $\mu\geq\kappa.$ Then a subset $A$ of $P_\kappa(\lambda)$ lies
in $NS_{\kappa,\lambda}^{[\mu]^{<\theta}}$ if and only if $B\cap\{a\in
C_g^{\kappa,\lambda}:a\cap\kappa\in\kappa\}=\phi$ for some $g : P_{\overline\theta.3}(\mu) \rightarrow
P_3(\lambda).$
}
\end{Lem}

The following is a straightforward generalization of a result of Foreman [1] :

\begin{Pro} {\sl Every $[\mu]^{<\theta}$-normal,
$(\mu^{<\overline\theta})^+$-saturated ideal on $P_\kappa(\lambda)$ is precipitous.
}
\end{Pro}

\section{$NS_{\kappa,\lambda}^{[\chi]^{<\theta}}\mid A = NS_{\kappa,\lambda}^{cf(\lambda)}\mid
A$}

\begin{Pro} 
\begin{enumerate}[i)]
\item {\sl Suppose $\lambda$  is a singular limit cardinal and $\theta$  a cardinal such
that $2\leq\theta\leq\kappa, \overline\theta\leq cf(\lambda)$ and
$\overline{cof}(NS_{\kappa,\tau}^{[\tau]^{<\theta}})\leq\lambda^{<\overline\theta}$ for every
cardinal $\tau$ with $\kappa\leq\tau<\lambda.$ Then there is
$A\in(NS_{\kappa,\lambda}^{[\lambda]^{<\theta}})^*$ such that
$NS_{\kappa,\lambda}^{[\lambda]^{<\theta}}=NS_{\kappa,\lambda}^{cf(\lambda)}\mid A.$
 }
\medskip
\item  {\sl Suppose $\lambda$ is a singular limit cardinal, $\theta$ is a cardinal such that
$2\leq\theta\leq\kappa,$ and $\chi$ is a cardinal such that
$\kappa\cdot(cf(\lambda))^+\leq\chi<\lambda$ and
$\overline{cof}(NS_{\kappa,\tau}^{[\chi]^{<\theta}})\leq\lambda^{<\overline\theta}$ for every
cardinal $\tau$ with $\chi\leq\tau<\lambda.$ Then there is
$A\in(NS_{\kappa,\lambda}^{[\lambda]^{<\theta}})^*$ such that
$NS_{\kappa,\lambda}^{[\chi]^{<\theta}}\mid  A = NS_{\kappa,\lambda}^{cf(\lambda)}\mid A.$
}
\end{enumerate}
\end{Pro}

{\bf Proof.} We prove both assertions simultaneously. Let us thus assume that $\lambda$ is a
singular limit cardinal, $\theta$ is a cardinal such that $2\leq\theta\leq\kappa,$ and $\chi$
 is a cardinal such that $\kappa\cdot(cf(\lambda))^+\leq\chi\leq\lambda$ and
$\overline{cof}(NS_{\kappa,\tau}^{[\chi\cap\tau]^{<\theta}})\leq\lambda^{<\overline\theta}$ for
every cardinal $\tau$ with $\pi\leq\tau<\lambda,$ where $\pi$ equals $\kappa$ if $\chi=\lambda,$
and $\chi$ otherwise. Let us also assume that $\overline\theta\leq cf(\lambda)$ in case
$\chi=\lambda.$ We are looking for $A\in(NS_{\kappa,\lambda}^{[\lambda]^{<\theta}})^*$ such that
$(NS_{\kappa,\lambda}^{[\chi]^{<\theta}}\mid A) = NS_{\kappa,\lambda}^{cf(\lambda)}\mid A.$
\\

Set $\mu=cf(\lambda)$ and select an increasing sequence of cardinals $<\lambda_\eta:\eta<\mu>$ so
that (a) $\displaystyle\bigcup_{\eta<\mu}\lambda_\eta=\lambda,$ (b) $\lambda_0\geq\kappa\cdot\mu,$
\ and (c) $\lambda_0\geq\chi$ in case $\chi<\lambda.$ For $\eta<\mu,$ pick a family $G_\eta$ of
functions from $P_{\overline\theta\cdot 3}(\chi\cap\lambda_\eta)$ to $P_3(\lambda_\eta)$ so that
$\mid G_\eta\mid\leq\lambda^{<\overline\theta}$ and for every
$H\in(NS_{\kappa,\lambda}^{[\chi\cap\lambda_\eta]^{<\theta}})^*,$ there is $y\in
P_\kappa(G_\eta)\setminus\{\phi\}$ such that $\displaystyle\{a\in\bigcap_{g\in
y}C_g^{\kappa,\lambda_\eta}:a\cap\kappa\in \kappa\}\subseteq H.$ Let
$\displaystyle\bigcup_{\eta<\mu}G_\eta=\{g_e:e\in P_{\overline\theta\cdot 3}(\lambda)\}.$
\\
Let $A$ be the set of all $a\in P_\kappa(\lambda)$ such that 
\begin{itemize}
\item $\overline\theta\subseteq a$ in case $\overline\theta<\kappa$ ;
\item $\omega\subseteq a$ ;
\item $a\cap\kappa\in\kappa$ ;
\item $k(\alpha)\in a$ for every $\alpha\in a,$ where $k : \lambda\rightarrow\mu$ is defined by
$k(\alpha) =$ the least $\eta<\mu$ such that $\alpha\in\lambda_\eta$ ;
\item If $\chi=\lambda,$ then $i(v)\in a$ for every $v\in
P_{|a\cap(\overline\theta\cdot 3)|}(a),$ where $i : P_{\overline\theta\cdot 3}(\lambda)\rightarrow\mu$ is
defined by $i(v) =$ the least $\eta<\mu$ such that $v\subseteq\lambda_\eta$ ;
\item $g_e(u)\subseteq a$ whenever $e\in P_{|a\cap(\overline\theta\cdot 3)|}(a)$ and $u\in
P_{|a\cap(\overline\theta\cdot 3)|}(a)\cap \hbox{dom}(g_e).$
\end{itemize}

It is immediate that $A\in (NS_{\kappa,\lambda}^{[\lambda]^{<\theta}})^*.$ Let us check that $A$
is as desired. Thus fix $B\in (NS_{\kappa,\lambda}^\mu)^+\cap P(A)$ and $f :
P_{\overline\theta\cdot 3}(\chi)\rightarrow P_3(\lambda).$ We must show that $B\cap
C_f^{\kappa,\lambda}\not=\phi.$ Given $\eta<\mu,$ define $p_\eta :
P_{\overline\theta\cdot 3}(\chi\cap\lambda_\eta) \rightarrow P_2(\lambda_\eta)$ by $p_\eta(v) =$ the least
$\sigma$ such that $\eta\leq\sigma<\mu$ and $f(v)\subseteq\lambda_\sigma.$ Also define $q_\eta :
P_{\overline\theta\cdot 3}(\chi\cap\lambda_\eta)\rightarrow P_3(\lambda_\eta)$ by $q_\eta(v) = 
\lambda_\eta\cap f(v).$ Select $x_\eta,y_\eta\in
P_\kappa(P_{\overline\theta\cdot 3}(\lambda))\setminus
\{\phi\}$ so that $\{g_e:e\in x_\eta\cup y_\eta\}\subseteq G_\eta, \{a\in
\displaystyle\bigcap_{e\in x_\eta}C_{g_e}^{\kappa,\lambda_\eta}:a\cap\kappa\in\kappa\}\subseteq
C_{p_\eta}^{\kappa,\lambda_\eta}$ and $\displaystyle\{a\in\bigcap_{e\in y_\eta}
C_{g_e}^{\kappa,\lambda_\eta}:a\cap\kappa\in\kappa\}\subseteq C_{q_\eta}^{\kappa,\lambda_n}.$ Now
pick $a\in B$ so that for any $\eta\in a\cap\mu, (\alpha) \  e\subseteq a$ for every $e\in
x_\eta\cup  y_\eta,$ and $(\beta)$ if $\overline\theta=\kappa,$ then $|e|<a\cap\kappa$ for every
$e\in x_\eta\cup y_\eta.$ Fix $v\in P_{|a\cap(\overline\theta\cdot 3)|}(a\cap\chi).$ There must be
$\eta\in a\cap\mu$ such that $v\subseteq\lambda_\eta.$ Then $a\cap\lambda_\eta\in
C_{p_\eta}^{\kappa,\lambda_\eta}$ since $x_\eta\subseteq P_{|a\cap(\overline\theta\cdot 3)|}(a).$
It follows that $v\cup f(v)\subseteq\lambda_\sigma$ for some $\sigma\in a\cap\mu.$ Now
$a\cap\lambda_\sigma\in C_{q_\sigma}^{\kappa,\lambda_\sigma},$ since $y_\sigma\subseteq
P_{|a\cap(\overline\theta\cdot 3)|}(a),$ so $f(v)\subseteq a.$ \hfill $\square$

In Proposition 2.1 (i) we assumed that $\overline\theta\leq cf(\lambda).$ Some condition of this
kind is necessary. In fact if  $u(\kappa,\lambda^{<\overline\theta})=\lambda^{<\overline\theta},$
then for each $A\in (NS_{\kappa,\lambda}^{[\lambda]^{<\theta}})^*,
NS_{\kappa,\lambda}^{[\lambda]^{<\theta}}\not= NS_{\kappa,\lambda}^{cf(\lambda)}\mid A$ since by
results of [4],

$$\overline{cof}(NS_{\kappa,\lambda}^{[\lambda]^{<\overline\theta}})
>\lambda^{<\overline\theta}\geq \lambda\geq\overline{cof}(NS_{\kappa,\lambda}^{cf(\lambda)}\mid
A).$$
\\

The following is immediate from Proposition 2.1 and Lemma 1.6 :
\\

\begin{Cor} {\sl Suppose $\chi>\kappa$ is a cardinal such that
$\overline{cof}(NS_{\kappa,\chi})\leq \chi^{+\kappa}.$  Then\break
$NS_{\kappa,\chi^{+\kappa}}^\chi\mid
 A = NS_{\kappa,\chi^{+\kappa}}^\kappa\mid A$ for  some $A\in(NS_{\kappa,\chi^{+\kappa}})^*.$
}
\end{Cor}

\section{Precipitousness}

\begin{Pro} {\sl Suppose $\lambda$ is a singular limit cardinal such that
$cf(\lambda)\geq\kappa$ and $\tau^{cf(\lambda)}<\lambda$ for every cardinal $\tau<\lambda.$ Then
there exists $B\in NS_{\kappa,\lambda}^*$ such that $NS_{\kappa,\lambda}^{cf(\lambda)}\mid B$ is
nowhere precipitous.
 }
 \end{Pro}

Proposition 3.1 will be obtained as a consequence of Lemmas 1.1 and 3.3.
\\

Let $\lambda$ be a singular limit cardinal of cofinality greater than or equal to $\kappa.$ Set
$\mu = cf(\lambda).$ Select a continuous, increasing sequence $<\lambda_\beta:\beta<\mu>$ of
cardinals so that (a) $\displaystyle\bigcup_{\beta<\mu}\lambda_\beta=\lambda,$ (b) $\lambda_0>\mu,$
and (c) $\lambda_0>2^\mu$ in case $\lambda>2^\mu.$ Le $E$ be the set of all limit ordinals
$\alpha<\mu$ with  $cf(\alpha)<\kappa.$ For $\alpha\in E,$ put $W_\alpha = \{a\in
P_\kappa(\lambda) : \cup a= \lambda_\alpha\}.$ Note that
$W_\alpha\in I_{\kappa,\lambda_\alpha}^*.$
\\

Let $B$ be the set of all $a\in P_\kappa(\lambda)$ such that (i) $0\in a,$ (ii) $\gamma+1\in a$
for every $\gamma\in a,$ (iii) $a\cap\kappa\in \kappa,$ (iv) $a\cap\mu = \{\beta\in\mu
:\lambda_\beta\in a\},$ and (v) for every $\gamma\in a,$ there is $\beta\in a\cap\mu$ such that
$\gamma<\lambda_\beta.$ Then clearly, $B\in NS_{\kappa,\lambda}^*.$ Moreover, $a\in
W_{\cup(a\cap\mu)}$ for every $a\in B.$ Note that for each $\alpha\in E, B\cap W_\alpha\subseteq
\{a\in P_\kappa(\lambda_\alpha):\cup(a\cap\mu) =\alpha\},$ so $B\cap W_\alpha\in
I_{\kappa,\lambda_\alpha}.$
\\

\begin{Lem}{\sl Suppose $u(\mu^+,\tau)<\lambda$ for every cardinal $\tau$ with
$\mu<\tau<\lambda.$ Then $\{\alpha\in E :|R\cap W_\alpha|\geq u(\mu^+,\lambda_\alpha)\}\in
NS_\mu^+$ for every $R\in (NS_{\kappa,\lambda}^\mu)^+\cap P(B).$
}
\end{Lem}

{\bf Proof.} Let us first show that for every $S\in (NS_{\kappa,\lambda}^\mu)^+\cap P(B),$ there is
$\alpha\in E$ such that $|S\cap W_\alpha|\geq u(\mu^+,\lambda_\alpha).$ Thus fix such an $S.$
Assume to the contrary that $|S\cap W_\alpha|<u(\mu^+,\lambda_\alpha)$ for every $\alpha\in E.$ For
$\alpha\in E,$ select $Z_\alpha\in I_{\mu^+,\lambda_\alpha}^+$ with $|Z_\alpha|<\lambda.$ Pick a
bijection $i : \displaystyle\bigcup_{\alpha<\mu}Z_\alpha\rightarrow\lambda$ and let $j$ denote the
inverse of $i.$ For $\alpha\in E,$ define $k_\alpha:P_\kappa(\lambda_\alpha)\rightarrow
P_{\mu^+}(\lambda_\alpha)$ by $k_\alpha(a) = \displaystyle\bigcup_{\beta\in a}(\lambda_\alpha\cap
j(\beta)),$ and select $y_\alpha\in P_{\mu^+}(\lambda_\alpha)$ so that $y_\alpha\setminus
k_\alpha(a)\not=\phi$ for every $a\in S\cap W_\alpha.$ Set $y = \displaystyle\bigcup_{\alpha\in
E} y_\alpha.$ Note that $y\in P_{\mu^+}(\lambda).$ For $\eta\in\mu,$ pick $z_\eta\in Z_\eta$ so
that $y\cap\lambda_\eta\subseteq z_\eta.$ Now let $D$ be the set of all $a\in P_\kappa(\lambda)$
such that $i(z_\eta)\in a$ for every $\eta\in a\cap\mu.$ Since $D\in(NS_{\kappa,\lambda}^\mu)^*,$
we can find $a\in S\cap D.$ Set $\alpha=\cup(a\cap\mu).$ Then $a\in W_\alpha$ and
$$y_\alpha\subseteq y\cap\lambda_\alpha =\bigcup_{\eta\in
a\cap\mu}(y\cap\lambda_\eta)\subseteq\bigcup_{\eta\in a\cap\mu}z_\eta=\bigcup_{\eta\in
a\cap\mu}j(i(z_\eta))\subseteq k_\alpha(a).$$
\\

Contradiction.
\\

It is now easy to show that the conclusion of the lemma holds. Fix $R\in
(NS_{\kappa,\lambda}^\mu)^+\cap P(B)$ and $T\in NS_\mu^*.$ Set $Q = \{a\in P_\kappa(\lambda)
:\cup(a\cap\mu)\in T\}.$ Since $Q\in (NS_{\kappa,\lambda}^\mu)^*,$ there must be some $\alpha\in
E$ such that $|(R\cap Q)\cap W_\alpha|\geq u(\mu^+,\lambda_\alpha).$ Then clearly, $\alpha\in T$
and $|R\cap W_\alpha|\geq u(\mu^+,\lambda_\alpha).$ \hfill$\square$
\\

\begin{Lem} {\sl Suppose $\tau^\mu<\lambda$ for every cardinal $\tau<\lambda.$ Then II has a winning strategy
in the game $G(NS_{\kappa,\lambda}^\mu|B).$
}
\end{Lem}

{\bf Proof.} For $g : P_3(\mu)\rightarrow P_2(\lambda)$ and $\alpha<\mu,$ define $g_\alpha:P_3(\mu)\rightarrow
P_2(\lambda_\alpha)$ by $g_\alpha(e) = \lambda_\alpha\cap g(e).$
\\

{\bf Claim 1.} {\sl Let $g : P_3(\mu)\rightarrow P_2(\lambda).$ Then $\{\alpha<\mu : B\cap W_\alpha\cap
C_{g_\alpha}^{\kappa,\lambda_\alpha}\subseteq C_g^{\kappa,\lambda}\}\in(NS_\mu\mid E)^*.$
}
\\

{\bf Proof of Claim 1.} Define $h : P_3(\mu)\rightarrow\mu$ by $h(e)=$ the least $\beta<\mu$ such that
$g(e)\subseteq\lambda_\beta.$ Let $Q$ be the set of all $\delta\in\mu$ such that $h(e)<\delta$ for
every $e\in P_3(\delta).$ Then clearly $Q\in NS_\mu^*.$ Now fix $\alpha\in E\cap Q$ and $a\in
B\cap W_\alpha\cap C_g^{\kappa,\lambda_\alpha}.$ Let $e\in P_3(a\cap\mu).$ Then $e\in
P_3(\alpha),$ so $g(e)\subseteq\lambda_\alpha.$ It follows that $g(e)\subseteq a,$ since
$\lambda_\alpha\cap g(e)\subseteq a.$ Thus $a\in C_g^{\kappa,\lambda}.$ \hfill$\square$
\\

{\bf Claim 2.} {\sl Let $X\in (NS_{\kappa,\lambda}^\mu)^+\cap P(B)$ and $Y\subseteq B.$ Suppose
that
$Y\cap W_\alpha\cap C_k^{\kappa,\lambda_\alpha}\not=\phi$ whenever $\alpha\in E$ and
$k:P_3(\mu)\rightarrow P_2(\lambda_\alpha)$ are such that $|X\cap W_\alpha\cap
C_k^{\kappa,\lambda_\alpha}|=\lambda_\alpha.$ Then $Y\in(NS_{\kappa,\lambda}^\mu)^+.$
}
\\

{\bf Proof of Claim 2.} Fix $g:P_3(\mu)\rightarrow P_2(\lambda).$ By Lemma 3.2 and Claim 1, there must be
$\alpha\in E$ such that $|(X\cap C_g^{\kappa,\lambda})\cap W_\alpha|=\lambda_\alpha^\mu$ and
$B\cap W_\alpha\cap C_{g_\alpha}^{\kappa,\lambda_\alpha}\subseteq C_g^{\kappa,\lambda}.$ Then
$Y\cap W_\alpha\cap C_{g_\alpha}^{\kappa,\lambda_\alpha}\not=\phi$ since $X\cap W_\alpha\cap
C_g^{\kappa,\lambda}\subseteq X\cap W_\alpha\cap C_{g_\alpha}^{\kappa,\lambda_\alpha}.$ Hence
$Y\cap C_g^{\kappa,\lambda}\not=\phi.$ \hfill$\square$
\\

For $\alpha\in E,$ consider the following two-person game $G_\alpha$ consisting of $\omega$ moves,
with player I making the first move : I and II alternately pick subsets of $B\cap W_\alpha,$ thus
building a sequence $<X_n:n<\omega>$ subject to the following two conditions : (1) $X_0\supseteq
X_1\supseteq\ldots,$ and (2) $X_{2n+1}\cap C_k^{\kappa,\lambda_\alpha}\not=\phi$ for every
$k:P_3(\mu)\rightarrow P_2(\lambda_\alpha)$ such that $|X_{2n}\cap
C_k^{\kappa,\lambda_\alpha}|=\lambda_\alpha^\mu.$ II wins the game if and only if
$\displaystyle\bigcap_{n<\omega}X_n=\phi.$
\\

{\bf Claim 3.} {\sl Let $\alpha\in E.$ Then II has a winning strategy $\tau_\alpha$ in the game
$G_\alpha.$
}
\\

{\bf Proof of Claim 3.} Let $Y_0, Y_1, Y_2,\ldots$ be the successive moves of player I. For
$n\in\omega,$ let $K_n$ be the set of all $k:P_3(\mu)\rightarrow P_2(\lambda_\alpha)$ such that $|Y_n\cap
C_k^{\kappa,\lambda_\alpha}|=\lambda_\alpha^\mu.$
\\

{\bf Case 1 :} $|K_n| =\lambda_\alpha^\mu$ for every $n<\omega.$ 

Given $n<\omega,$ set $K_n =
\{k_{n,\xi}:\xi<\lambda_\alpha^\mu\}$ and let $\tau_\alpha(Y_0,\ldots, Y_n)
=\{y_{n,\xi}:\xi<\lambda_\alpha^\mu\},$ where $y_{n,\xi}\in (Y_n\cap
C_{k_{n,\xi}}^{\kappa,\lambda_\alpha})\setminus\{y_{q,\zeta}:q<n$ and $\zeta\leq\xi\}.$
\hfill$\square$\\

{\bf Case 2 :} There is $n\in\omega$ such that $|K_n|<\lambda_\alpha^\mu.$

Let $m$ be the least such $n.$ Define $\tau_\alpha(Y_0,\ldots,Y_m)$ so that
$|\tau_\alpha(Y_0,\ldots,Y_m)|<\lambda_\alpha^\mu,$ and set
$\tau_\alpha(Y_0,\ldots,Y_m,Y_{m+1})=\phi.$

Finally, consider the strategy $\tau$ for player II in $G(NS_{\kappa,\lambda}^\mu|B)$ defined by
$\tau(X_0)=\displaystyle\bigcup_{\alpha\in E}\tau_\alpha(X_0\cap B\cap W_\alpha)$ and for $n>0,$
$$\tau(X_0,\ldots,X_n)=\bigcup_{\alpha\in E}\tau_\alpha(X_0\cap B\cap W_\alpha, X_1\cap
W_\alpha,\ldots,X_n\cap W_\alpha).$$ Using Claims 2 and 3, it is easy to check that the strategy II
is a winning one.\hfill $\square$

\begin{Pro} {\sl Let $\lambda$ be a singular limit cardinal with
$cf(\lambda)\geq\kappa,$ and $\theta$ be a cardinal with $2\leq\theta\leq\kappa.$ Suppose that for
every cardinal $\tau$ with $\kappa\leq\tau<\lambda, \tau^{cf(\lambda)}<\lambda$ and
$\overline{cof}(NS_{\kappa,\tau}^{[\tau]^{<\theta}})\leq\lambda^{<\overline\theta}.$ Then the
ideal $NS_{\kappa,\lambda}^{[\lambda]^{<\theta}}$ is nowhere precipitous.
}
\end{Pro}

{\bf Proof.} By Propositions 2.1 and 3.1, one can find $A, B\in
(NS_{\kappa,\lambda}^{[\lambda]^{<\theta}})^*$ such that
$NS_{\kappa,\lambda}^{[\lambda]^{<\theta}}=NS_{\kappa,\lambda}^{cf(\lambda)}|A$ and
$NS_{\kappa,\lambda}^{cf(\lambda)}|B$ is nowhere precipitous. Now for every
$T\in(NS_{\kappa,\lambda}^{[\lambda]^{<\theta}})^+,$
$$NS_{\kappa,\lambda}^{[\lambda]^{<\theta}}|T =
(NS_{\kappa,\lambda}^{[\lambda]^{<\theta}}|B)|T=((NS_{\kappa,\lambda}^{cf(\lambda)}|A)|B)|T =
(NS_{\kappa,\lambda}^{cf(\lambda}|B)|(A\cap T),$$
where $A\cap T\in(NS_{\kappa,\lambda}^{cf(\lambda)}|B)^+.$\hfill$\square$

Universit\'e de Caen - CNRS \\
Laboratoire de Math\'ematiques \\
BP 5186 \\
14032 CAEN CEDEX\\
France\\
matet@math.unicaen.fr\\

Institute of Mathematics \\
The Hebrew University of Jerusalem\\
91904 Jerusalem\\
Israel \\

and\\

Department of Mathematics\\
Rutgers University\\
New Brunswick, NJ 08854\\
USA\\
shelah@math.huji.ac.il

\end{document}